\documentclass[12pt]{article}
\usepackage{amsmath,amssymb}
\newcommand{\G}{\Gamma}
\newcommand{\nin}{\noindent}
\newcommand{\vs}{\vspace*}

\newcommand{\eset}{\emptyset}

\newcommand{\seq}{\subseteq}

\newcommand{\llrarrow}{\longleftrightarrow}
\newcommand{\lrarrow}{\leftrightarrow}

\newcommand{\ol}{\overline}

\begin{document}

\baselineskip=16pt
\begin{center}
{\bf\LARGE Sub-semigroups determined by the zero-divisor graph }

\vskip 20pt

Wu Tongsuo {\footnote {Corresponding author, tswu@sjtu.edu.cn}} \\

Department of Mathematics\\
 Shanghai Jiaotong  University\\
   Shanghai 200030, P.R.China\\

   and

    Lu Dancheng  {\footnote { ludancheng@suda.edu.cn}} \\

 Department of Mathematics\\
 Suzhou University\\
 Suzhou 215006, P. R. China\\

\vskip 30pt

\begin{minipage}{12cm}

\vs{3mm}{\bf Abstract. }{\small In this paper we study
sub-semigroups of a zero-divisor semigroup $S$ determined by
properties of the zero-divisor graph $\G(S)$. We use these
sub-semigroups to study the correspondence between zero-divisor
semigroups and zero-divisor graphs. We study properties of
sub-semigroups of Boolean semigroups via the zero-divisor graph. As
an application, we provide a characterization of the graphs which
are zero-divisor graphs of Boolean rings.}

\vs{3mm}\nin {\small Key Words:} {\small zero-divisor semigroup,
sub-semigroup, zero-divisor graph, Boolean graph }
\end{minipage}
\end{center}

\vs{4mm}\begin{center}{\bf 1. Introduction}\end{center}

\vs{3mm} For any commutative semigroup $S$ with zero element $0$,
there is an undirected zero-divisor graph $\G(S)$ associated with
$S$. The vertex set of $\G(S)$ is the set of all nonzero
zero-divisors of $S$,  and for distinct vertices $x$ and $y$ of
$\G(S)$, there is an edge linking $x$ and $y$ if and only if $xy=0$
(\cite{FR}). In \cite{FR} and \cite{FL}, some fundamental properties
and possible structures of $\G(S)$ were studied. For example, for
any semigroup $S$, it was proved that $\G(S)$ is a connected simple
graph with diameter less than or equal to 3, that the core of
$\G(S)$ is a union of triangles and squares while any vertex of
$\G(S)$ is either an end vertex or in the core, {\it if there exists
a cycle in $\G(S)$}, and that for any non-adjacent vertices $x,y$,
there exists a vertex $z$ such that $N(x)\cup N(y)\seq \ol{N(z)}$.
In \cite{FL}, the authors also provided a descending chain of ideals
$I_k$ of $S$, where $I_k$ consists of all elements of $S$ with
vertex degree greater than or equal to $k$ in $\G(S)$. These general
structure results are very helpful when one studies the following
general problem: Given a connected simple graph $G$, does there
exist a semigroup $S$ such that $\G(S)\cong G$? Some classes of
graphs were given in \cite{FL,FR},\cite{WCH,WULU,ZWU} to give
positive or negative answers to this problem. The zero-divisor graph
was also extensively studied for commutative rings, see, e.g.,
\cite{AL,DF,B,AN,M}.

In this paper, we study some sub-semigroup structures of a
zero-divisor semigroup by properties of the graph $\G(S)$. First we
prove that if there is an vertex $x$ such that the graph has three
parts $\G(S)=T_x\cup\{x\}\cup C_x$, where $T_x$ contains all end
vertices adjacent to $x$ such that there is no edge linking a vertex
in $T_x$ with a vertex in $C_x$, and that either $C_x$ is non-empty
or $x$ is adjacent to at least one end vertex,  then $\{x\}\cup
C_x\cup \{0\}$ is a sub-semigroup of $S$. In particular, if $x$ is
adjacent to at least one end vertex and $G$ has a cycle, then
$x^2=0$ or $x^2=x$. We also provide various conditions such that the
end vertices adjacent to a single vertex together with 0 forms a
sub-semigroup. For Boolean semigroups (i.e., semigroups with
$x^2=x,\forall x\in S$), we found a descending chain of
sub-semigroups which in some sense is dual to the $I_k$ mentioned
above. As an application, we characterize the graphs which are
zero-divisor graphs of Boolean rings. We also use these
sub-semigroups to study the correspondence between semigroups and
zero-divisor graphs. In particular, we construct a kind of graph $G$
which has a unique corresponding zero-divisor semigroup, where the
core of $G$ is the union of a square and a triangle. Two vertices of
$G$ can be adjacent to arbitrarily many (finite or infinite) end
vertices, while if one adds end vertices to the other three
vertices, the resulting graph has no corresponding semigroups.

All semigroups in this paper are multiplicatively commutative {\it
zero-divisor} semigroups with zero element $0$, where $0x=0$ for all
$x\in S$, and all graphs in this paper are undirected simple and
connected. For a given connected simple graph $G$, if there exists a
zero-divisor semigroup $S$ such that $\G(S)\cong G$, then we say
that {\it $G$ has corresponding semigroups}, and we call $S$ a {\it
semigroup determined by the graph $G$}. For any vertices $x,y$ in a
graph $G$, if $x$ and $y$ are adjacent, we denote it as $x-y$ or
occasionally, $x\lrarrow y$. A vertex is called an {\it end vertex},
if its vertex degree is one. The {\it core} of a graph $G$, which
will be denoted as $K(G)$, is the largest subgraph of $G$ every edge
of which is an edge of a cycle in $G$. We denote the complete graph
with $n$ end vertices by $K_n$. Similarly, we denote a complete
bipartite graph with two parts of sizes $m,n$ by $K_{m,n}$.

\vs{4mm}\begin{center} {\bf 2. Sub-semigroups related to a single
vertex}
\end{center}

\vs{3mm}\nin{\bf Theorem 2.1.} {\it Let S be a zero-divisor
semigroup with graph G. Assume that for an element x of S, there
exists a subset $T_x$ of $S-\{0,x\}$ satisfying the following
conditions}:

(1) {\it $T_x$ contains all end vertices adjacent to x}.

(2) {\it There is no edge linking a vertex in $T_x$ with a vertex in
$S-(T_x\cup \{0,x\})$  when $S-(T_x\cup \{0,x\})\not=\eset$}.

(3) {\it Either $S-(T_x\cup \{0,x\})\not=\eset$, or $G$ has a cycle
and x is adjacent to at least one end vertex}.

\nin{\it Then $T=S-T_x$ is a sub-semigroup of $S$,}

\vs{3mm}\nin{\bf Proof.} First assume that $S-(T_x\cup
\{0,x\})=\eset$. Then by assumption (3), $x$ is adjacent to at least
one end vertex $y\in T_x$ and $x$ is in the core of $G$. By
\cite[Theorem 4]{FL}, $x^2$ is not an end vertex and in particular
$x^2\not=y$. In this case, $T=S-T_x=\{x,0\}$ must be a sub-semigroup
of $S$. In fact, if $x^2\not=0$ and $x^2\not=x$, then we have a path
$x^2-y-x$, contradicting to the fact that $y$ is an end vertex.

Now let us assume that $S-(T_x\cup \{0,x\})$ is not empty. Then
there is an element $x\not=z\in T$ such that $z-x$. First let us
consider $xT$. If for some $t\in T$, $xt=y\in T_x$, then we have a
path $z-y-x$ where $z\not\in T_x, x\not=z, y\in T_x$, contradicting
with the condition (2). Thus we must have $xT\subseteq T$. For any
$a\in T-\{x\}$, we now proceed to prove $aT\subseteq T$. Assume in
the contrary that there is an element $x\not=b\in T$ such that
$ab=y\in T_x$. If $a$ is not adjacent to $x$, then there is an
element $c\in T-\{a,x\}$ such that $a-c$. In this case there is a
path $c-y-x$, where $y\in T_x, c\not\in T_x$, contradicting again
with the condition (2). If $a$ is adjacent to $x$, then by
assumption (1), $a$ is not an end vertex. Thus by condition (2)
there is an element $c\in T-\{a,x\}$ such that $a-c$. In this case,
we also have a path $c-y-x$, where $y\in T_x, c\not\in T_x$. In each
case, there is a contradiction with the condition (2). These
contradictions show that $aT\subseteq T$. Finally, we obtain
$T^2\subseteq T$ and hence $T$ is a sub-semigroup of $S$. \quad
$\Box$

\vs{3mm}\nin{\bf Corollary 2.2.} {\it Let S be a zero-divisor
semigroup with graph G. For a vertex x of G, assume that $T_x$ is
nonempty and denote $T=S-T_x$, where $T_x$ is the end vertices
adjacent to x. Then $T$ is a non-trivial sub-semigroup of $S$. If in
addition $G$ has a cycle, then $\{x,0\}$ is a sub-semigroup of $S$.}

\vs{3mm}\nin{\bf Proof.} When $T_x$ is not empty, we have $x\in
T\subset S$. The set $T_x$ of all end vertices adjacent to x
certainly satisfies all conditions of Theorem 2.1. Thus the
corollary follows from Theorem 2.1.\quad $\Box$

\vs{3mm}For a vertex $v$ of a connected graph $G$, if $v$ is not
an end vertex and there is no end vertex adjacent to $v$, then $v$
is said to be an {\it internal } vertex of $G$. For an internal
vertex $v$ of a graph $G$ and another graph $H$ which is disjoint
with $G$, we get a new graph {\it by attaching the graph $H$ to
$v$, } i.e., we add an edge linking $v$ to every vertex of $H$.

 \vs{3mm}\nin{\bf
Corollary 2.3.} {\it If G is a graph which is not the graph of any
semigroup and $F$ is a graph obtained from G by attaching another
graph to an internal vertex of G, then F is not the graph of any
semigroup. In particular, if G is a graph which is not the graph
of any semigroup and $F$ is a graph obtained from G by adding some
ends to internal vertices of G, then F is not the graph of any
semigroup.}

\vs{3mm}\nin{\bf Proof.} This follows immediately from Theorem 2.1.
\quad $\Box$

\vs{3mm}\nin{\bf Example 2.4.} (\cite[Corollary 4.2(1)]{WCH})  For
any $n\ge 4$, let $M_{n,k}$ be the complete graph
$K_n=\{a_1,a_2,\cdots, a_n\}$ together with $k$ end vertices
$\{x_i|1\le i\le k\}$ such that $a_i-x_i$ ($n\ge k\ge 4$). Then
$M_{n,k}$ has no corresponding semigroups.

\vs{3mm}\nin{\bf Proof.} By \cite[Theorem 2.2]{WULU}, the graph
$M_{n,3}$ has no corresponding semigroups. Since $M_{n,k}$ is
obtained by adding some end vertices to $M_{n,3}$, the result
follows from Corollary 2.3.\quad $\Box$

\vs{3mm}Theorem 2.1 and Corollary 2.2 are particularly powerful
when the related graph has a unique corresponding zero-divisor
semigroup, as is illustrated by the following examples.

\vs{3mm}\nin{\bf Example 2.5.} By \cite[Corollary 2.4]{WCH}, the
following graph $G$ has a unique corresponding zero-divisor
semigroup (The multiplication table is written on the right hand.):
$$\begin{array}{ccccccc}
\circ&\llrarrow &\circ&\llrarrow&\circ\\
\Big\updownarrow&&\Big\updownarrow&\diagup\\
\circ&\llrarrow&\circ\\
&&\tiny{\text{ (G)}}&\\
\end{array},\quad
\begin{array}{c|cccccc}
\cdot& a_1&a_2&a_3&x_1&x_2\\
\hline
a_1&0&0&0&0&a_1\\
a_2&&0&0&a_2&0\\
a_3&&&0&a_2&a_1\\
x_1&&&&x_1&0\\
x_2&&&&&x_2\\
&&&\tiny{\text{ Tab. 1}}&&
\end{array}.
$$
Now Let us consider graphs $H$ obtained by adding some end vertices
to some vertices of $G$. We have the following conclusions:

(1) Each of the following graphs has a unique zero-divisor
semigroup:
\[\begin{array}{ccccccccccccccc}
&&$U$&&&\\
&&\Big\updownarrow&&\\
\circ&\llrarrow &\circ&\llrarrow&\circ\\
\Big\updownarrow&&\Big\updownarrow&\diagup&\\
\circ&\llrarrow&\circ&\llrarrow&V\\
&&\tiny{\text{Fig.1}}&&&
\end{array}.\]
\nin where both $U$ and $V$ consists of some end vertices.

(2) For other graphs $H$ obtained by adding some end vertices to
some vertices of $G$, $H$ has no corresponding semigroups.

\vs{3mm}\nin{\bf Proof.} (i) Consider the following graphs:
\[\begin{array}{cccccccccccccccccc}
\circ&\llrarrow &\circ&\llrarrow&\circ&&&&&&\circ&\llrarrow&\circ&\llrarrow\circ\\
\Big\updownarrow&&\Big\updownarrow&\diagup&\Big\updownarrow&&&&&&\Big\updownarrow&&\Big\updownarrow&\diagup\\
\circ&\llrarrow&\circ&&U&,&&&U&\llrarrow&\circ&\llrarrow&\circ&&\\
&&\tiny{\text{Fig.2}}&&&&&&&&&\tiny{\text{Fig.3}}&&
\end{array},\]
where $U$ is nonempty and it consists of some end vertices. In
each case, there exist non-adjacent vertices $x,y$ such that for
any $z$, $N(x)\cup N(y)\not\subseteq \ol{N(z)}$. By \cite [Theorem
1(1)]{FL} and Theorem 2.1, only graphs of Fig. 1 may have
semigroups.

(ii) Now consider the graphs $H$ in Fig.1 and let us first assume
$U=\emptyset$, and $a_2u=0,\forall u\in V$. For any $u\in V$, from
$$x_2u\in ann(x_1)\cap ann(a_2)\subseteq \{x_2,a_1\}$$
\nin we have $x_2u=x_2$, since $x_2u=a_1$ implies
$0=a_1a_3=(a_3x_2)u=a_1u$  by Corollary 2.2. Hence $u^2\not=0$,
and $a_1=a_1x_2=a_1u$. Since
$$x_1u\in ann(x_1)\cap ann(a_1)\cap ann(a_2)=\{a_2\},$$
we must have $x_1u=a_2$. This implies $x_1u^2=0$. Since
$u^2\not=0,x_1^2=x_1$, we have either $u^2=a_1$ or $u^2=x_2$. But
$u^2=a_1$ implies that $a_1=a_1x_2=(ux_2)u=x_2u=x_2$, a
contradiction. Thus we must have $u^2=x_2,\forall u\in V$. Now
that $x_1u=a_2$ implies $x_1(a_3u)=0$, thus we must have
$a_3u=a_1$, since otherwise one has $a_3u=x_2$, and this will
imply $a_1x_2=0$. Finally, for any $u\not=v\in V$, we conclude
$uv=x_2$. In fact, since $uv\in ann(a_2)$, $uv\not=x_1$. If
$uv=a_1$, then $a_1=a_1u=u^2v=x_2v=x_2$. If $uv=a_2 $, then
$0=a_2u=u^2v=x_2$. If $uv=a_3$, then $0=(a_1u)v=a_1$. If $uv=u$,
then $0=ua_2=u(vx_1)=ux_1=a_2$. These contradictions show that
$uv=x_2$.

By Corollary 2.2, The graph $H$ in Fig.1 uniquely determines the
multiplication table on $S=V(H)\cup \{0\}$. The final work is to
verify that the table defines an associative binary operation on
$S$. This is really the case. In the following Tab.1, we list the
table for $U=\emptyset, V=\{u,v\}$ with $|V|=2$.
$$\begin{array}{c|ccccccc}
\cdot& a_1&a_2&a_3&x_1&x_2&u&v\\
\hline
a_1&0&0&0&0&a_1&a_1&a_1\\
a_2&&0&0&a_2&0&0&0\\
a_3&&&0&a_2&a_1&a_1&a_1\\
x_1&&&&x_1&0&a_2&a_2\\
x_2&&&&&x_2&x_2&x_2\\
u&&&&&&x_2&x_2 \\
v&&&&&&&x_2\\
&&&&&\tiny{\text{Tab.2}}&&
\end{array}.
$$

(iii) Now assume $U=\{u_i\,|\,1\le i\le r\}$,  $V=\{v_j\,|\,1\le
j\le s\}$ and $a_1u=0,\forall u\in U,a_2v=0,\forall v\in V$. By
Corollary 2.2, we need only to determine the value of $uv$ for any
$u\in U, v\in V$. Since
$$uv\in ann(a_1)\cap ann(a_2)\subseteq \{a_1,a_2,a_3\},$$ \nin we
must have $uv=a_3$. In fact, if $uv=a_1$, then
$0=x_1a_1=(x_1v)u=a_2u$, a contradiction. Finally, it is routine
to test the associative law. Below we list the table for
$U=\{u_1,u_2\}, V=\{v_1,v_2\}$ and this will end the proof:
$$\begin{array}{c|ccccccccccccccc}
\cdot& a_1&a_2&a_3&x_1&x_2&u_1&u_2&v_1&v_2\\
\hline
a_1&0&0&0&0&a_1&0&0&a_1&a_1\\
a_2&&0&0&a_2&0&a_2&a_2&0&0\\
a_3&&&0&a_2&a_1&a_2&a_2&a_1&a_1\\
x_1&&&&x_1&0&x_1&x_1&a_2&a_2\\
x_2&&&&&x_2&a_1&a_1&x_2&x_2\\
u_1&&&&&&x_1&x_1&a_3&a_3\\
u_2&&&&&&&x_1&a_3&a_3 \\
v_1&&&&&&&&x_2&x_2\\
v_2&&&&&&&&&x_2\\
&&&&&\tiny{\text{Tab.3}}&&&&
\end{array}.
$$
\quad $\Box$

\vs{3mm}\nin{\bf Example 2.6.} The following graph has
corresponding semigroups:
\[\begin{array}{cccccccc}
U&\llrarrow &\circ\\
&&\Big\updownarrow&\diagdown\\
V&\llrarrow&\circ&\llrarrow&\circ&\llrarrow&W\\
&&&\text{{\tiny Fig.4}} &
\end{array}\]
\nin where $U,V$ and $W$ consist of end vertices.

\vs{3mm}\nin{\bf Proof.} For $U=\{x_1,u\}, V=\{x_2,v\},
W=\{x_3,w\}$, the following associative table defines a semigroup
whose zero-divisor graph is isomorphic to the graph of Fig.4:
$$\begin{array}{c|cccccccccccccccccc}
\cdot& a_1&a_2&a_3&x_1&u&x_2&v&x_3&w\\
\hline
a_1&a_1&0&0&0&0&a_1&a_1&a_1&a_1\\
a_2&&a_2&0&a_2&a_2&0&0&a_2&a_2\\
a_3&&&a_3&a_3&a_3&a_3&a_3&0&0\\
x_1&&&&x_1&x_1&a_3&a_3&a_2&a_2\\
u &&&&&x_1&a_3&a_3&a_2&a_2\\
x_2&&&&&&x_2&x_2&a_1&a_1\\
v&&&&&&&x_2&a_1&a_1\\
x_3&&&&&&&&x_3&x_3\\
w&&&&&&&&&x_3\\
&&&&&\tiny{\text{Tab.4}}&&&&
\end{array}.
$$
\quad $\Box$

\vs{3mm}\nin{\bf Proposition 2.7.} {\it Let $G$ be the graph of a
semigroup S and assume that G has a cycle. For a vertex x of the
graph G which is not an end vertex, let $T_x$ be the end vertices
adjacent to $x$. Then $T_x\cup \{0\}$ is a sub-semigroup of $S$, if
$x^2\not=0$. }

\vs{3mm}\nin{\bf Proof.} For any $u,v\in T_x$ with $uv\not= 0$,
obviously $uv\not=x,uv-x$ since $x^2\not=0.$ If $uv=c$ is not an end
vertex of $G$, then it is in the core of $G$ since the graph $G$ has
a cycle. Then there exists another vertex $d\not=x,c$ in the core
such that $x-c-d$. Then from $0=dc=(du)v$, we obtain $du=x$ by
\cite[Theorem 4]{FL}. But then we have $x^2=d(xu)=0$, a
contradiction. Thus $uv$ must be an end vertex and hence $uv\in
T_x$. \quad $\Box$

\vs{3mm}We remark that the condition of $x^2\not=0$ could not be
dropped from Proposition 2.7.

\vs{3mm}\nin{\bf Example 2.8.} Consider the semigroup
$S=\{0,a_1,a_2,a_3,x_1\}$ with the multiplication table
\begin{center}\vs{2mm}$\begin{array}{c|ccccccccccccccccccccccccc}
\cdot& a_1&a_2&a_3&x_1 \\
\hline
a_1&0&0&0&0             \\
a_2&0&a_1&0&a_1             \\
a_3&0&0&a_3&a_3&             \\
x_1&0&a_1&a_3&a_3&              \\
&&&\text{\tiny Tab.5}&& \\
\end{array}$
\end{center}
Obviously $T_{a_1}\cup \{0\}=\{0,x_1\}$ is not a sub-semigroup of
$S$. \quad$\Box$

\vs{3mm}\nin{\bf Theorem 2.9.} {\it Let $G$ be the graph of a
semigroup S. Let $s,t$ be two vertices of $G$ such that
$T_s\not=\eset,T_t\not=\eset$ and that $s^2=0,t^2=0 $, where $T_s$
is the end vertices adjacent to $s$. Then the following conclusions
hold:}

(1) {\it $T_sT_t\seq N(s)\cap N(t)$ and $sS=\{0,s\}, tS=\{0,t\}$.}

(2) {\it If one of the following conditions holds, then $T_s\cup
\{0\}$ is a reduced sub-semigroup of $S$:} \mbox(I) {\it $u^2=0,
\forall u\in N(s)\cap K(G)$, where $K(G)$ is the core of $G$.}
\mbox(II) {\it The edge $s-t$ is not contained in any
quadrilaterals.}

\vs{3mm}\nin{\bf Proof.} (1) Let us first show that $T_sT_t\seq
N(s)\cap N(t)$. Assume that $y\in T_s, x\in T_t$, and consider the
following
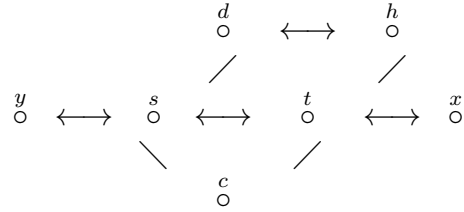
\[\begin{array}{cccccccccccc}
&&&\stackrel{d}{\circ}&\llrarrow &\stackrel{h}{\circ}&&\\
&&&\diagup&&\diagup\\
\stackrel{y}{\circ}&\llrarrow &\stackrel{s}{\circ}&\llrarrow&\stackrel{t}{\circ}&\llrarrow&\stackrel{x}{\circ}\\
&&\diagdown&&\diagup &&&&&\\
&&&\stackrel{c}{\circ}&&&
\end{array}\]
Since $xy\in \mbox{ann}(s)\cap \mbox{ann}(t)$, The value of $xy$
has only three possibilities, i.e., $xy=s,xy=t, xy=c \in N(s)\cap
N(t)$.

If $xy=s$, then we have $0\not=sx=x^2y$, and hence $x^2\not=0$.
Since $0=ts=(ty)x$, we have $ty=t$. Thus $ty^2=t$ and hence
$y^2\not=0$. On the other hand, we have $0=sy=xy^2$ and thus
$y^2=t$. Finally, we have $t^2=ty^2=ty=t$, contradicting with the
assumption. Thus $xy\not=s$. By symmetry, we also have $xy\not=t$.
Thus we must have $xy=c\in N(s)\cap N(t)$. In this case we have
$x^2\not=0$ since $0\not=cx=x^2y$. Similarly we also have
$y^2\not=0$. Therefore $sS=$\mbox{ann}$(y)=\{0,s\}$, and
$tS=\mbox{ann}(x)=\{0,t\}$. In particular, we have $sx=s$ and
$ty=t$.

(2) Now we show that $T_s\cup \{0\}$ is a reduced sub-semigroup of
$S$. For any $y\in T_s$, we already know that $y^2\not=0$ by case
3. Now assume $y^2=d\in \ol{N(s)}$. Then we have $td=ty^2=t$, thus
$td^2=t$ and hence $d^2\not=0$. Hence $d\not=s,t$. Also $d\not=c
$, since $c-t$ but $d$ is not adjacent to $t$. If condition (I)
holds, then we know that $d=y^2\in T_s$ is an end vertex. If
condition (II) holds, then $d$ must be an end vertex adjacent to
$s$, since otherwise we have another vertex $h\in K(G)$ such that
$h\not= s$ and $d-h$. In this case, we have a square $s-d-h-t-s$,
a contradiction to the assumption (II).

 For distinct elements $y,z\in
T_s$, obviously $yz=d\in \ol{N(s)}$. If $d=s$, then $0=sz=z^2y$.
Then either $z^2=y$ or $z^2-y$, contradictions. Thus $yz\not=s$
always holds. Now assume $yz=d\not=s$. First assume that the
condition (I) holds. Since $0=(sy)z$, we must have $sy=s,\forall
y\in T_s$. If $d$ is not an end vertex, then $0=d^2=(dy)z$. Thus
$dy=s$, which implies $dy^2=sy=0$. Then we have $d-y^2-s$,
contradicting to the fact that $y^2$ is an end vertex. In the
other case we assume that the condition (I) holds. In this case,
we also have $dt=t$. Thus if $d\not\in T_s$, then there exists an
$h\in G(K)$ such that $h\not=s$ and $hd=0$. Then we will obtain a
square $s-d-h-t-s$, contradicting with the assumption (II) again.
Finally, in each case, we have $yz\in T_s$. Thus $T_s\cup \{0\}$
is a reduced sub-semigroup of $S$. \quad$\Box$

\vs{3mm} For any $1\le m\le |V(G)|$, the graph $G$ is said to be
{\it $m$-uniquely determined (by neighborhoods)}, if for any $x,y\in
V(G)$ with $|N(x)|=|N(y)|=m$, $N(x)=N(y)$ implies $x=y$.

\vs{3mm}\nin{\bf Proposition 2.10.} {\it Let $G$ be the graph of a
finite zero-divisor semigroup S. Let $s$ be a vertex of $G$ with
the greatest degree $m$. Assume that $G$ is $m$-uniquely
determined and $s^2=s$. Then $\{s,0\}$ is an ideal of $S$. }

\vs{3mm}\nin{\bf Proof.} For any $x\in S$ with $xs\not=0$, we have
$sx\not\in N(s)$. Thus $N(s)\seq N(sx)$ and hence $N(s)=N(sx)$. By
assumption, we obtain $sx=s$. This shows that $Ss\seq \{s,0\}$.
\quad $\Box$

\vs{4mm}\begin{center} {\bf 3. Sub-semigroups of Boolean
semigroups\\
and uniquely determined graphs}
\end{center}

\vs{3mm} Recall that a ring $R$ (a semigroup $R$) is called a {\it
Boolean ring }({\it a Boolean semigroup}), if $r^2=r,\forall r\in
R$. It is easy to see that vertices of the zero-divisor graph
$\G(R)$ of a Boolean ring $R$ are uniquely determined by their
neighborhoods, i.e., $N(x)=N(y)$ only if $x=y$ for any vertices of
$\G(R)$ (\cite [Page 2]{LUWU}), where $N(x)=\{y\in \G(R)\,|\,y-x\}$
is the neighbors of $x$ in $\G(R)$. However, this result does not
hold for general Boolean semigroups.

\vs{3mm}\nin{\bf Proposition 3.1.} (\cite[Propsosition 3.2]{ZWU})
{\it There exists an Boolean semigroup $S$ such that $\G(S)$ is a
complete $r$-partite graph. The vertices of this $\G(S)$ are not
uniquely determined by their neighborhoods when $r>1$ and at least
one part of $\G(R)$ contains more than one vertices. }

\vs{3mm}\nin{\bf Proof.} In fact, Let $S=\{0\}\cup (\cup_{i=1}^r
A_i)$ be a disjoint union of $r+1$ nonempty subsets, where
$A_i=\{a_{ik_i}\,|\,1\le k_i\le m_i\}$. Define
$a_{ik_i}^2=a_{ik_i}$, $a_{ir}a_{is}=a_{i1}$ ($r\not=s$),
$a_{ik}a_{jl}=0$ ($i\not=j$). Then $S$ is a commutative zero-divisor
semigroup whose zero-divisor graph $\G(S)$ is the complete
$r$-partite graph. Obviously, $S$ is a Boolean semigroup.\quad
$\Box$

\vs{3mm} For any zero-divisor semigroup $S$, define an equivalent
relation $\sim$ in $S^*$ in the following natural way:

\centerline{ $ x\sim y$ if and only if $N(x)=N(y)$.}

\nin For any $0\not=x\in S$, denote by $S_x=\{y\in S\,|\, y\not=0,
N(y)=N(x)\}$ the equivalent class containing $x$. Assume $S^*$ has
$m_e$ classes.

\vs{3mm}\nin{\bf Theorem 3.2.} {\it Let $S$ be a Boolean semigroup
with zero-element 0}.

 (1) {\it  Then $S_x$ is a sub-semigroup of S, and $0\not\in
S_x$. }

(2) {\it Let $S_{\leqslant x}=\{y\in S^*\,|\,N(y)\seq N(x)\}$.
Then $S_{\leqslant x}$ is a sub-semigroup of $S$, $0\not\in
S_{\leqslant x}$ and $S_x$ is an ideal of $S_{\leqslant x}$.}

\vs{3mm}\nin{\bf Proof.} (1) For any $z,y\not=x$, if
$N(y)=N(x)=N(z)$, then  obviously $xy\not=0$, and $N(x)\seq N(xy)$
since $xy\not\in N(x)$. Now for any $u\in N(xy)$, we have
$u\not=x,y,xy$. If $xu\not=0$, then $xu\not= y$. Therefore  $xu\in
N(y)=N(x)$, and hence $xu=x^2u=0$, a contradiction. Thus we must
have $N(xy)\seq N(x)$ and hence $N(x)=N(xy)$. This shows that $xy\in
S_x,\forall x,y\in S_x$. Since $N(z)=N(x)=N(xy)=N(y)$, we have
$N(x)=N(xyz)=N(yz)$. Thus $yz\in S_x$. Thus $S_x$ is a sub-semigroup
of $S$.

(2) For any $x\not=y\in S_{\leqslant x}$, if $xy=0$, then $x\in
N(y)$ but $x\not\in N(x)$, contradicting with the assumption of
$N(y)\seq N(x)$. Thus $xy\not=0$. Since $x^2=x$, we have $xy\not\in
N(x)$, and therefore $N(x)\seq N(xy)$. Conversely, for any $u\in
N(xy)$, we have $uxy=0, u\not=x,y,xy$. If $ux\not=0$, then
$ux\not=y$. Thus $ux\in N(y)\seq N(x)$, and hence $ux=ux^2=0$, a
contradiction. Thus we also have $N(xy)\seq N(x)$. This shows that
$N(x)=N(xy)$. Now if $z\in S_{\leqslant x}$, then we have $N(yz)\seq
N(xyz)=N(x)$ and hence $yz\in S_{\leqslant x}$. This shows that
$S_{\leqslant x}$ is a sub-semigroup of $S$, $0\not\in S_{\leqslant
x}$ and $S_x$ is an ideal of $S_{\leqslant x}$. \quad $\Box$

\vs{3mm}\nin{\bf Corollary 3.3.} {\it For any Boolean semigroup $S$
with zero-element 0, the zero-divisor graph $\G(S)$ is uniquely
determined if and only if $N(y)\seq N(x)$ implies $yx=x, \forall
y,x\in S^*$. }

\vs{3mm}\nin{\bf Proof.} $\Longleftarrow.$ If $N(x)=N(y)$, then
$x=xy=y$.

$\Longrightarrow.$ If $N(y)\seq N(x)$, then by the proof of Theorem
3.2, we have $N(xy)=N(x)$. If $\G(S)$ is uniquely determined, then
we obtain $xy=x$. \quad $\Box$

\vs{3mm} As an application of Theorem 3.2 and Corollary 3.3, we
obtain the following results which should be compared with Example
3.1:

\vs{3mm}\nin{\bf Corollary 3.4.} {\it For any $n\ge 3$, each
connected subgraph $G$ of the complete graph $K_{n+1}$ containing
the complete graph $K_n$ has a unique zero-divisor Boolean semigroup
if $|G|\not=|K_n|+(n-1)$. }

\vs{3mm}\nin{\bf Proof.} Assume that there exists a Boolean
semigroup $S$ such that $\G(S)= G$. Assume $S-\{0\}=V(G)$. We only
prove the following case. Assume that $G$ is the complete graph
$K_n=\{a_i\,|\,1\le i\le n\}$ together with an end vertex $x_1$,
where $a_1-x_1$. Since $N(x_1)\seq N(a_i)$ ($i\not=1$), by Theorem
3.2(2) we obtain $N(a_i)=N(a_ix_1)$. Since $G$ is uniquely
determined by neighborhoods, we have $a_i=a_ix_1$. Thus we have a
unique associative multiplication table on the vertices of $G$ such
that $\G(S)=G$. Notice that $G$ is uniquely determined by
neighborhoods if and only if the degree of $x_1$ is not $n-1$. When
the degree of $x_1$ is $n-1$, the graph $G$ has multiple
corresponding Boolean semigroups such that $\G(S)\cong G$. \quad
$\Box$

\vs{3mm} Recall that a {\it two-star graph} is a graph consists of
two star graphs with exactly one edge connecting the two centers. By
\cite[Theorem 1.3]{FR}, the two-star graphs and their connected
subgraphs are all the possible zero-divisor graphs of semigroups
that contains no cycle.

\vs{3mm}\nin{\bf Corollary 3.5.} {\it The following graphs have no
Boolean semigroups:}

(1) {\it The complete graph $K_n$ together with more than one end
vertices ($n\ge 4$). }

(2) {\it The complete bipartite graph $K_{m,n}$ together with end
vertices ($m,n\ge 2$)}.

(3) {\it Any two star graph which is not a star graph.}

\vs{3mm}\nin{\bf Proof.} Assume that there exists a Boolean
semigroup $S$ such that $\G(S)= G$.

(1) and (2). Now assume $G$ is the complete bipartite graph
$$K_{m,n}=\{a_i\,|\,1\le i\le m\}\cup \{b_j\,|\,1\le j\le n\}$$
\nin together with an end vertex $x_1$, where $a_1-x_1$ ($m,n\ge
2$). Obviously, $a_1a_2=a_1$. Now consider $a_2x_1$. If
$a_2x_1=a_1$, then $a_1^2=0$. Thus $a_2x_1=a_k$ for some $k\not=1$.
But then we have $(a_1a_2)x_1\not=a_1(a_2x_1)$. The negative results
of other cases follows from \cite[Theorem 3.2]{WULU} and
\cite[Corollary 3.4]{WCH}.

(3) Assume that $G=\{x_i\,|\,m\ge i\ge 1\}\cup \{a,b\}\cup
\{y_j\,|\,1\le j\le n\}$, where $x_i-a-b-y_j, \forall i,j$. If
$S=V(G)\cup\{0\}$  is a Boolean semigroup such that $\G(S)=G$, then
by Theorem 3.2(2) we have $x_1b=b,y_1a=a$ since $N(x_1)\seq N(b),
N(y_1)\seq N(a)$. Then we have $0=by_1=(x_1y_1)b$ and
$0=ax_1=(x_1y_1)a$. Thus we obtain
$$x_1y_1\in \{a, y_j\}\cap\{b,x_i\},$$
\nin which is obviously impossible. \quad $\Box$

\vs{3mm}Recall that a semigroup $S$ is called a {\it reduced}
semigroup, if there exist no nonzero nilpotent elements in $S$.

\vs{3mm}\nin{\bf Proposition 3.6.} {\it For a reduced zero-divisor
semigroup $S$, if the graph $\G(S)$ is uniquely determined, then $S$
is a Boolean semigroup.}

\vs{3mm}\nin{\bf Proof.} For any $0\not=x\in S$, obviously $N(x)\seq
N(x^2)$. Conversely, if $u\in N(x^2)$, then we have $(ux)^2=0$. Thus
$u\not=x$ and $ux=0$. Thus $N(x^2)\seq N(x)$. \quad $\Box$

\vs{3mm} Recall that a Boolean algebra is {\it a distributive
complemented} lattice with the smallest element $0$ and the largest
element $1$. Recall from \cite[Page 225]{DF} that a graph $G$ is
called {\it complemented}, if for each vertex $x$ of $G$, there
exists a vertex $y$ such that $x\perp y$, i.e., $x\not=y$ and the
edge $x-y$ is not part of any triangle. $G$ is called {\it uniquely
complemented}, if $G$ is complemented and whenever $x\perp y$ and
$x\perp z$, then $N(y)=N(z)$. Following \cite{LUWU}, a simple graph
$G$ is called a {\it Boolean graph}, if $G$ is the zero-divisor
graph of some Boolean ring. By \cite[Theorem 4.2, Corollary
4.3]{LUWU}, a Boolean graph with more than two vertices has a unique
corresponding zero-divisor semigroup and a unique corresponding
ring. In the final part of this paper, we use Theorem 3.2 to give a
characterization of Boolean graphs.

\vs{3mm}\nin{\bf Theorem 3.7.} {\it A simple connected graph $G$ is
a Boolean graph if and only if the following conditions hold:}

(1) {\it The graph $G$ is uniquely determined.}

(2) {\it The graph $G$ is uniquely complemented.}

(3) {\it For any $x,y\in V(G)$ with $N(x)\cap N(y)\neq \emptyset$,
there exists some $z\in V(G)$ such that $N(x)\cap N(y)=N(z)$.}

(4) {\it The graph $G$ has a corresponding Boolean semigroup.}
\vs{3mm}

\noindent {\bf Proof.} $\Longrightarrow $. It is obvious.

$\Longleftarrow $. By (4), we can assume that $V(G)\cup \{0\}$ is a
Boolean semigroup. Denote $\emptyset$ and $V(G)$ by  $N(1)$ and
$N(0)$ respectively. Set $R=\{0,1\}\cup V(G)$ and let $P(G)=
\{N(x)|x\in R \}$.  Then $P(G)$ has a natural partial order.  First
we show that $N(x)\vee N(y)=N(xy)$. Let $z\in V(G)\cup \{0\}$ such
that $N(x)\cup N(y)\subseteq N(z)\subseteq N(xy)$. If $xy\not=0$,
then $z\not=0$. By the previous Corollary 3.3, we have $z=zx=zy$ and
so $z=zxy=xy$, $N(z)=N(xy)$. If $xy=0$ and $z\not=0$, then we have
$z=zx=zy$ and thus $z=z(xy)=0$, a contradiction. Thus in each case,
we have $N(x)\vee N(y)=N(xy)$. Second we define $N(x)\wedge
N(y)=N(x)\cap N(y)$. By condition (3), $P(G)$ becomes a lattice. We
claim that $P(G)$ is a distributive lattice. In fact, it suffice to
prove that
$$(N(x)\wedge N(y))\vee N(z)=(N(x)\vee N(z))\wedge (N(y)\vee
N(z)).$$ Of course, we have $(N(x)\wedge N(y))\vee N(z)\subseteq
(N(x)\vee N(z))\wedge (N(y)\vee N(z)).$ Conversely, let $t \in
(N(x)\vee N(z))\wedge (N(y)\vee N(z))$ and assume $N(x)\wedge
N(y)=N(k)$. Then $txz=0$ and $tyz=0$. So $tz\in N(x)\wedge N(y)$ and
it follows that $tzk=0$. This implies that $t\in N(k)\vee N(z)$,
proving the claim. By condition (3), we already know that $P(G)$ is
a Boolean algebra. In the following we denote by $N(\overline{x})$
the complement of $N(x)$.

Clearly, there is a multiplication operation on $R$. Define an
addition $+$ on $R$ as follows: Given $x,y\in R$, then there is a
unique $z\in R$ such that $$(N(x)\vee N(\overline{y}))\wedge
(N(\overline{x})\vee N(y))=N(z).$$ \nin We then define $x+y=z$. It
is routine to verify that $R$ is a Boolean ring and that $\G(R)=G$.
This completes the proof. \quad$\Box$

\vs{3mm}


\begin{thebibliography}{gg}

\bibitem{AL}
D.F. Anderson and P.S. Livingston,  The zero-divisor graph of a
commutative ring, {\it  J. Algebra}, {\bf 217}(1999), 434-447.

\bibitem{DF}
D.F. Anderson, Ron Levy,  Jay Shapiro, Zero-divisor graphs, von
Neumann regular rings, and Boolean algebras, {\it J. Pure Applied
Algebra} {\bf 180}(2003), 221-241.

\bibitem{B}
I. Beck, Coloring of commutative rings, {\it J. Algebra}, {\bf
116}(1988),  208-226.

\bibitem{AN}
D.D. Anderson and M. Naseer, Beck's coloring of a commutative ring,
{\it J. Algebra} {\bf 159}(1993), 500-514.

\bibitem{FL}
F. DeMeyer, L. DeMeyer, Zero-divisor graphs of semigroups, {\it J.
Algebra} {\bf 283}(2005)190-198.

\bibitem{FR}
F.R. DeMeyer, T. McKenzie, and K. Schneider, The zero-divisor graph
of a commutative semigroup, {\it Semigroup Forum} {\bf 65}(2002),
206-214.

\bibitem{LUWU}D.C. Lu and T.S. Wu, The zero-divisor graphs which are
uniquely determined by neighborhoods. Preprint 2005.

\bibitem{M}
B. Mulay, Cycles and symmetries of zero-divisors, {\it Comm.
Algebra} {\bf 30}(7)(2002), 3533-3558.

\bibitem{WCH}
T.S. Wu and L. Chen, Simple graphs and commutative zero-divisor
semigroup, {\it Comm. Algebra}.

\bibitem{WULU}
T.S. Wu and D.C. Lu, Zero-divisor semigroups and some simple graphs,
{\it Comm. Algebra} (To appear).

\bibitem{ZWU}
M. Zuo and T.S. Wu, A new graph structure of commutative semigroups,
{\it Semigroup Forum} 70:1(2005), 71---80.  DOI:
10.1007/s00233-004-0139-8

\end{thebibliography}
\end{document}